\documentclass[ijoc,nonblindrev]{informs3}
\usepackage{mathrsfs}
\usepackage{latexsym}
\usepackage{amsmath}
\usepackage{theorem}
\usepackage{amssymb}
\usepackage{multirow}
\usepackage{graphicx}
\usepackage{enumerate}
\usepackage{url}
\usepackage{cases}
\usepackage{lscape}
\usepackage{mathtools}
\usepackage{algorithmic}
\usepackage{algorithm}
\OneAndAHalfSpacedXII 

\TheoremsNumberedThrough     

\EquationsNumberedThrough    

\newenvironment{preuve}{\noindent\textit{Proof. \ }}{\hfill$\Box$\medskip\par }
\newcounter{triang}

\newcounter{binary}

\newcommand{\ipopt}{\texttt{Ipopt}}
\newcommand{\mosek}{\texttt{Mosek}}

\newcommand{\baron}{\texttt{Baron}}

\newcommand{\bb}{branch-and-bound}
\newcommand{\sbb}{spatial branch-and-bound}

\newcommand{\rcopf}{\texttt{RC-OPF}}
\newcommand{\opfdiag}{\texttt{COPF}}

\newcommand{\diag}{\mathrm{d}}
\newcommand{\Diag}{\mathrm{D}}

\newcommand{\setMC}{\cal{D}}

\newcommand{\setN}{\cal{N}}
\newcommand{\setC}{\cal{C}}
\newcommand{\setE}{\cal{E}}

\newcommand{\OPT}{(P)}

\newcommand{\OPF}{(OPF)}

\makeatletter
\let\markeq\ltx@label
\makeatother

\usepackage{url}

\begin{document}


\TITLE{A tight compact quadratically constrained  convex relaxation of the Optimal Power Flow problem}

\RUNTITLE{A tight compact quadratically constrained  convex relaxation of the OPF problem}
\ARTICLEAUTHORS{\AUTHOR{Am\'elie Lambert} \AFF{Cnam-CEDRIC, 292 Rue St Martin FR-75141 Paris Cedex 03, \EMAIL{amelie.lambert@cnam.fr}
}}

\RUNAUTHOR{Am\'elie Lambert}


\ABSTRACT{In this paper, we consider  the optimal power flow (OPF) problem which consists in determining the power production at each bus of an electric network by minimizing the production cost.  Our contribution is an exact solution algorithm for the OPF problem. It consits in a spatial branch-and-bound algorithm based on  a compact quadratically constrained convex relaxation. This compact relaxation is computed  by solving the rank relaxation once at the beginning of the algorithm. The key point of this approach is that the lower bound at the root node of the branch-and-bound tree is equal to the rank relaxation value, but is obtained by solving a quadratic convex problem which is much faster than solving a SDP. To construct this compact relaxation, we add only  $\mathcal{O}(n)$  variables that model the squares of the initial variables, where $n$ is the number of buses in the power system. The relations between the initial and auxiliary variables are therefore non-convex. By relaxing them in our relaxation, we have only $\mathcal{O}(n)$  equalities to force by the branch-and-bound algorithm to prove global optimality. Our first experiments show that our new algorithm \texttt{Compact OPF} (\opfdiag) performs better than the methods of the literature we compare it with. 
}


\maketitle

\section{Introduction}

\par The Optimal Power Flow (OPF) problem consists in the determination of the power production at different buses of an electric network that minimizes a production cost. The electrical transmission network is modeled by a graph $G = (\setN,\setE)$. Each network point belongs to the set $\setN$ of nodes (i.e the set of buses), and their connections (i.e. the set of transmission lines) are modeled by the set of edges $\setE$. We assume that there is an electric demand at each node also called load. We distinguish two classes of nodes: $\setN =\setN_g \cup \setN_d$, where  $\setN_g$ is the set of nodes that generates and flows the power (the generator nodes), and $\setN_d$ is the set of nodes that only flows the power (the consuming nodes). The aim of the OPF problem is to satisfy demand of all buses while minimizing the total production costs of the generators such that the solution obeys Ohm's law and and Kirchhoff's law.

\par This problem is naturally formulated with complex variables. Let $ \mathbf{Y} \in  \mathbb{C}^{n\times n}$, where  $ \lvert \setN \rvert = n$ be the admittance matrix, which has component $ \mathbf{Y}_{i,k} =  \mathbf{R}_{i,k} + \textrm{j} \mathbf{I}_{i,k}$ for each line $(i,k)$ of the network, and  $\mathbf{R}_{i,i} = \mathbf{ r}_{i,i} - \displaystyle{\sum_{i \neq k}} \mathbf{R}_{i,k}$, $\mathbf{I}_{i,i} =  \mathbf{i}_{i,i} - \displaystyle{\sum_{i \neq k}} \mathbf{I}_{i,k}$, where $\mathbf{r}_{i,i}$ (resp. $\mathbf{i}_{i,i}$) is the shunt conductance (resp. susceptance) at bus $i$, and $\textrm{j}^2 = -1$. Let $p_i, q_i$ be the real and reactive power output of the generator node $i$, and $\mathbf{p}_i, \mathbf{q}_i$ the given real and reactive power output of the load node $i$. We consider here the complex voltage in the rectangular form: $V_i = e_i + \textrm{j}f_i$ where $\lvert V_i \rvert^2 = e^2_i + f^2_i$ is the voltage magnitude, and we denote by $\delta(i)$ the set of adjacent nodes of bus $i$. With the above notation, the OPF problem can be modeled by the well known rectangular formulation of~\cite{Torres1998AnIM}:

\begin{small}
\begin{numcases}{(OPF)} 
\min  h(p) =  \displaystyle{\sum_{i \in \setN_g}} \Big ( \mathbf{C}_{i,i} p^{2}_i  + \mathbf{c}_ip_i \Big ) \nonumber \\
\mbox{s.t. } \nonumber \\
\quad p_i - \mathbf{p}_i =   \mathbf{R}_{i,i} (e_i^2 + f_i^2)  +  \displaystyle{\sum_{k \in \delta(i)}} \Big [\mathbf{R}_{i,k}(e_ie_k + f_if_k) - \mathbf{I}_{i,k}(e_if_k - e_kf_i) \Big ]&$i \in \setN_g$ \label{real_eq_g} \\
\quad - \mathbf{p}_i =   \mathbf{R}_{i,i} (e_i^2 + f_i^2)  +  \displaystyle{\sum_{k \in \delta(i)}} \Big [\mathbf{R}_{i,k}(e_ie_k + f_if_k) -   \mathbf{I}_{i,k}(e_if_k - e_kf_i) \Big ]&$i \in \setN_d$ \label{real_eq_d} \\
\quad q_i - \mathbf{q}_i =   -\mathbf{I}_{i,i} (e_i^2 + f_i^2)  +  \displaystyle{\sum_{k \in \delta(i)}} \Big [ - \mathbf{I}_{i,k}(e_ie_k +f_if_k) -   \mathbf{R}_{i,k}(e_if_k - e_kf_i) \Big ]&$i \in \setN_g$ \label{ima_eq_g} \\
\quad - \mathbf{q}_i =   -\mathbf{I}_{i,i} (e_i^2 + f_i^2)  +  \displaystyle{\sum_{k \in \delta(i)}} \Big [ - \mathbf{I}_{i,k}(e_ie_k + f_if_k) -   \mathbf{R}_{i,k}(e_if_k - e_kf_i) \Big ]&$i \in \setN_d$ \label{ima_eq_d} \\
\quad \underline{\mathbf{v}}_i \leq e_i^2 + f_i^2 \leq \overline{\mathbf{v}}_i & $i \in \setN$ \label{volt_magn} \\
\quad   \underline{\mathbf{p}}_i \leq p_i \leq  \overline{\mathbf{p}}_i & $i \in \setN_g$ \label{real_bound}\\
   \quad   \underline{\mathbf{q}}_i \leq q_i \leq  \overline{\mathbf{q}}_i & $i \in \setN_g$ \label{ima_bound}\\
\quad(e,f) \in (\mathbb{R}^{n},\mathbb{R}^{n}), \, \, (p,q) \in (\mathbb{R}^{\lvert \setN_g \rvert},\mathbb{R}^{\lvert \setN_g \rvert}) &  \label{setR}
\end{numcases}
\end{small}
\medskip

\noindent where  $\mathbf{C} \in \mathbf{S}^+_{\lvert \setN_g \rvert}$ is a diagonal and semi-definite matrix, $\mathbf{c} \in \mathbb{R}^{\lvert \setN_g \rvert}$ is the vector of linear costs of the power injection at each generator node, $( \underline{\mathbf{v}}, \overline{\mathbf{v}})  \in( \mathbb{R}^{n} ,\mathbb{R}^{n})$ are the bounds on the voltage magnitude, and $ (\underline{\mathbf{p}}, \overline{\mathbf{p}},\underline{\mathbf{q}}, \overline{\mathbf{q}})  \in (\mathbb{R}^{\lvert \setN_g \rvert},\mathbb{R}^{\lvert \setN_g \rvert}, \mathbb{R}^{\lvert \setN_g \rvert},\mathbb{R}^{\lvert \setN_g \rvert})$. This formulation has $2(n +2\lvert \setN_g \rvert)$ variables, $2n$ quadratic equalities~(\ref{real_eq_g})-(\ref{ima_eq_d}) that enforces  the active and reactive power balances at each node, $2n$ quadratic inequalities~(\ref{volt_magn}) that models the voltage magnitude, and $4\lvert \setN_g \rvert$ box constraints~(\ref{real_bound})-(\ref{ima_bound}).

\medskip

The first results of the literature for solving the  OPF problem were focused on optimal local solutions, mostly by adapting interior point methods, see, e.g.,~\cite{Wu1993ADN,Torres1998AnIM,JCC02,WMZT07}. In the context of global optimization, one requires furthermore to determine lower bounds on the OPF problem. For this, the second-order cone programming (SOCP) and the semidefinite programming (SDP) relaxations were first used (see~\cite{Jab06,BWFW08, BW09, LavLow12}). The most used SDP relaxation, also named the rank relaxation, leads to very tight lower bounds on the OPF problem. In particular, it was proven in~\cite{GLTL15} that this relaxation is exact  for a restricted class of problems and under some assumptions. In other cases, where it is not exact, it can be strengthen until the optimal solution value, following the ideas of the hierarchy of moment relaxation problems~\cite{Las01,Par03}, that can be applied to any polynomial optimisation problem. This approach was specialized in the context of the OPF problem in~\cite{JMPG15}, and showed its efficiency to solve small-size problems. It was also used in~\cite{MolHis15} to strengthen the lower bounds for larger problems. Unfortunately, in practice, using interior point methods for solving several large SDP relaxations, which sizes increase at each rank of the hierarchy, is intractable for large networks. Several specialized algorithms that exploit the sparsity of power networks were thus proposed as in~\cite{Jab06,Jab12,MHLD13, MolHis15, MSL15}. 

More recently, several cheaper computable convex relaxations were introduced for the OPF problem. For instance, linear and quadratic envelopes for trigonometric functions in
the polar formulation of the OPF problem are constructed in~\cite{CofVan14,CHV15b,CHV17}, and strong SOCP relaxations were introduced in~\cite{KDS16}. These polynomial bounds can then be used within a spatial branch-and-bound framework to solve the problem to global optimality. 

Another exact solution approach called \rcopf~was proposed in~\cite{Godard19,EGLMR19}. This method is a specialization of the  Quadratic Convex Reformulation approach proposed in~\cite{EllLam19} for generic quadratic problems. It consits in a branch-and-bound algorithm based on the computation of a tight quadratic convex relaxation the OPF problem that is computed  by solving the rank relaxation once at the beginning of the algorithm. The key point of this approach is that the lower bound at the root node of the branch-and-bound tree is equal to the rank relaxation value, but is obtained by solving a quadratic convex problem which is much faster than solving a SDP. This tight quadratic convex relaxation has a quadratic convex objective function and linear constraints. Moreover, it has a size of $\mathcal{O}((2n)^2)$, since it relies on the introduction of $(2n)^2$ additional variables that model all the possible products of the original variables. Then a spatial branch-and-bound is performed which aim is to enforce the equality between each additional variable and its corresponding product. Unfortunately, and as illustrated in Section~\ref{result}, in practice enforcing these $(2n)^2$ equalities can be very time consuming for the OPF problem.
\medskip

\par Our contribution is an exact solution algorithm for the OPF problem that is based on a compact quadratic convex relaxation. This relaxation  has only $\mathcal{O}(n)$ auxiliary variables and constraints. Moreover,  as in method  \rcopf, its optimal value also equals to the optimal value of the well known rank relaxation. We thus get a cheaper computable convex relaxation than is method \rcopf, that is as tight as the rank relaxation. Finally, to solve $(OPF)$ to global optimality, we perform a spatial \bb~algorithm based on this new quadratic convex relaxation. A key advantage of our compact relaxation is that we only have to enforce the  equality between  $2n$ additional variables and their corresponding products to prove global optimality. Finally, we evaluate our method \texttt{Compact OPF} (\opfdiag) on several instances of the literature and we show that our algorithm outperforms method \rcopf~and the generic global optimization solver \baron~\cite{baron}.

\medskip
\par The paper is organized as follows. In Section~\ref{conv_cont}, we describe how we convexify the quadratic constraints of $(OPF)$. Then, in Section~\ref{family}, we present a new family of equivalent quadratic convex formulations to $(OPF)$. In Section~\ref{best_relax}, we show how we calculate the best quadratic convex relaxation within this family. Finally, in Section~\ref{result}, we present some computational results. Section~\ref{conc} draws a conclusion. 

\section{A compact convexification of the quadratic constraints}
\label{conv_cont}

\par We start by observing that the structure of the formulation $(OPF)$ is specific. First, only variables $p$ are involved into the objective function. Moreover, the matrix $\mathbf{C}$ is diagonal and positive semidefinite, hence the objective function is convex and separable. It follows that the non convexities only come from the quadratic constraints~(\ref{real_eq_g})-(\ref{volt_magn}). More precisely, in the constraints, variables $e$ and $f$ are only involved into quadratic forms, while variables $p$ and $q$ only in linear forms. 
\medskip

Starting from the latter observations, our idea is to build an equivalent problem to $(OPF)$, where the original constraints are convexified thanks to $2n$ auxiliary variables $z = (z^e,z^f) \in \mathbb{R}^{2n} $ that model the squares of the initial variables $e$ and $f$:
\begin{numcases}{}
 \quad z_i^e = e_i^2  & $i \in \setN$    \label{prod_e}\\
\quad  z_i^f = f_i^2  & $i \in \setN$   \label{prod_f}
\end{numcases}
Using these new variables it is easy to rewrite Constraints~(\ref{volt_magn}) into a convex form by simply linearizing them, and we get:
\begin{align}
  \underline{\mathbf{v}}_i \leq z_i^e + z_i^f \leq \overline{\mathbf{v}}_i   \quad i \in \setN \label{volt_magn_c} 
\end{align}
\medskip

We now handle the equality constraints. For this, the first step is to transform each equality~(\ref{real_eq_g})-(\ref{ima_eq_d}) into two inequalities. We thus introduce for all $i \in \setN_g$,  Constraints~(\ref{real_eq_g1})-(\ref{ima_eq_d2}) obviously equivalent to Constraints~(\ref{real_eq_g})-(\ref{ima_eq_d}):
\begin{numcases}{}
p_i - \mathbf{p}_i - \mathbf{R}_{i,i} (e_i^2 + f_i^2)  -  \displaystyle{\sum_{k \in \delta(i)}} \Big [\mathbf{R}_{i,k}(e_ie_k + f_if_k) -   \mathbf{I}_{i,k}(e_if_k - e_kf_i) \Big ] \leq 0  &$i \in \setN_g$ \label{real_eq_g1} \\
- p_i + \mathbf{p}_i +   \mathbf{R}_{i,i} (e_i^2 + f_i^2)  +  \displaystyle{\sum_{k \in \delta(i)}} \Big [\mathbf{R}_{i,k}(e_ie_k + f_if_k) -   \mathbf{I}_{i,k}(e_if_k - e_kf_i) \Big ] \leq 0&$i \in \setN_g$  \label{real_eq_g2} \\
- \mathbf{p}_i - \mathbf{R}_{i,i} (e_i^2 + f_i^2)  -  \displaystyle{\sum_{k \in \delta(i)}} \Big [\mathbf{R}_{i,k}(e_ie_k + f_if_k) -   \mathbf{I}_{i,k}(e_if_k - e_kf_i) \Big ] \leq 0 & $i \in \setN_d$ \label{real_eq_d1} \\
\mathbf{p}_i +   \mathbf{R}_{i,i} (e_i^2 + f_i^2)  +  \displaystyle{\sum_{k \in \delta(i)}} \Big [\mathbf{R}_{i,k}(e_ie_k + f_if_k) -   \mathbf{I}_{i,k}(e_if_k - e_kf_i) \Big ] \leq 0  & $i \in \setN_d$\label{real_eq_d2} \\
q_i - \mathbf{q}_i +   \mathbf{I}_{i,i} (e_i^2 + f_i^2)  +  \displaystyle{\sum_{k \in \delta(i)}} \Big [ \mathbf{I}_{i,k}(e_ie_k + f_if_k) +   \mathbf{R}_{i,k}(e_if_k - e_kf_i) \Big ] \leq 0   &$i \in \setN_g$ \label{ima_eq_g1} \\
- q_i + \mathbf{q}_i -   \mathbf{I}_{i,i} (e_i^2 + f_i^2)  -  \displaystyle{\sum_{k \in \delta(i)}} \Big [ \mathbf{I}_{i,k}(e_ie_k + f_if_k) +  \mathbf{R}_{i,k}(e_if_k - e_kf_i) \Big ] \leq 0& $i \in \setN_g$  \label{ima_eq_g2}  \\
- \mathbf{q}_i  +   \mathbf{I}_{i,i} (e_i^2 + f_i^2)  +  \displaystyle{\sum_{k \in \delta(i)}} \Big [  \mathbf{I}_{i,k}(e_ie_k + f_if_k) +   \mathbf{R}_{i,k}(e_if_k - e_kf_i) \Big ] \leq 0   & $i \in \setN_d$  \label{ima_eq_d1}\\
\mathbf{q}_i -   \mathbf{I}_{i,i} (e_i^2 + f_i^2)  -  \displaystyle{\sum_{k \in \delta(i)}} \Big [  \mathbf{I}_{i,k}(e_ie_k + f_if_k) +   \mathbf{R}_{i,k}(e_if_k - e_kf_i) \Big ] \leq 0  & $i \in \setN_d$ \label{ima_eq_d2}
\end{numcases}

Then, recall that Inequalities~(\ref{real_eq_g1})-(\ref{ima_eq_d2}) are only quadratic on variables $e$ and $f$. To make them convex, we apply the smallest eigenvalue method introduced in~\cite{HamRub70}.  Denote by $A^1_i$ (resp. $A^2_i$, $A^3_i$, $A^4_i$) the sub-matrix of the Hessian of the $i^{th}$ Constraints~(\ref{real_eq_g}) (resp.~(\ref{real_eq_d}), (\ref{ima_eq_g}), (\ref{ima_eq_d})) that corresponds to the quadratic terms involving variables $e$ and $f$ only. Let  $\lambda(A^1_i) $ be the smallest eigenvalue of matrix $A^1_i$, and $\diag (\lambda(A^1_i))$ be the diagonal matrix where each diagonal term equals $\lambda(A^1_i)$. To rewrite Inequality~(\ref{real_eq_g1}) (resp.~(\ref{real_eq_g2}))  as a convex function, we add the quadratic quantity $ - \lambda(A^1_i) \displaystyle{\sum_{k\in \setN}} (e_k^2 + f_k^2 - z_k^e -z_k^f)$ (resp. $- \lambda(-A^1_i) \displaystyle{\sum_{k\in \setN}}  (e_k^2 + f_k^2 - z_k^e -z_k^f)$) to it. As a consequence the Hessian matrix of the new function is $A^1_i -  \diag (\lambda(A^1_i) )$ (resp. $-A^1_i -  \diag (\lambda(-A^1_i) )$) that is obviously a positive semidefinite matrix. Moreover, the value of the convexified function remains the same as soon as for all $k \in \setN$, $e_k^2 + f_k^2 - z_k^e -z_k^f=0$, or equivalently when Equalities~(\ref{prod_e})-(\ref{prod_f}) are satisfied. We thus obtain the set of convex Constraints~(\ref{real_eq_g_c1})-(\ref{ima_eq_d_c2}):

\begin{scriptsize}
\begin{numcases}{}
p_i - \mathbf{p}_i - \mathbf{R}_{i,i} (e_i^2 + f_i^2)  -  \displaystyle{\sum_{k \in \delta(i)}} \Big [\mathbf{R}_{i,k}(e_ie_k + f_if_k) -   \mathbf{I}_{i,k}(e_if_k - e_kf_i) \Big ]-  \lambda(A^1_i)\displaystyle{\sum_{k\in \setN}} \Big  ( e_k^2 + f_k^2 - z^e_k - z^f_k \Big ) \leq 0  & \hspace{-0.4cm} $i \in \setN_g$ \label{real_eq_g_c1} \\
- p_i + \mathbf{p}_i +   \mathbf{R}_{i,i} (e_i^2 + f_i^2)  +  \displaystyle{\sum_{k \in \delta(i)}} \Big [\mathbf{R}_{i,k}(e_ie_k + f_if_k) -   \mathbf{I}_{i,k}(e_if_k - e_kf_i) \Big ]-  \lambda(-A^1_i)  \displaystyle{\sum_{k\in \setN}} \Big  ( e_k^2 + f_k^2 - z^e_k - z^f_k \Big ) \leq 0& \hspace{-0.4cm} $i \in \setN_g$  \label{real_eq_g_c2} \\
- \mathbf{p}_i - \mathbf{R}_{i,i} (e_i^2 + f_i^2)  -  \displaystyle{\sum_{k \in \delta(i)}} \Big [\mathbf{R}_{i,k}(e_ie_k + f_if_k) -   \mathbf{I}_{i,k}(e_if_k - e_kf_i) \Big ] -  \lambda(A^2_i)\displaystyle{\sum_{k\in \setN}} \Big  ( e_k^2 + f_k^2 - z^e_k - z^f_k \Big )\leq 0 &  \hspace{-0.4cm} $i \in \setN_d$ \label{real_eq_d_c1} \\
\mathbf{p}_i +   \mathbf{R}_{i,i} (e_i^2 + f_i^2)  +  \displaystyle{\sum_{k \in \delta(i)}} \Big [\mathbf{R}_{i,k}(e_ie_k + f_if_k) -   \mathbf{I}_{i,k}(e_if_k - e_kf_i) \Big ]  -  \lambda(-A^2_i)  \displaystyle{\sum_{k\in \setN}} \Big  ( e_k^2 + f_k^2 - z^e_k - z^f_k \Big )\leq 0  & \hspace{-0.4cm}  $i \in \setN_d$\label{real_eq_d_c2} \\
q_i - \mathbf{q}_i +   \mathbf{I}_{i,i} (e_i^2 + f_i^2)  +  \displaystyle{\sum_{k \in \delta(i)}} \Big [ \mathbf{I}_{i,k}(e_ie_k + f_if_k) +   \mathbf{R}_{i,k}(e_if_k - e_kf_i) \Big ] -  \lambda(A^3_i)  \displaystyle{\sum_{k\in \setN}}\Big  ( e_k^2 + f_k^2 - z^e_k - z^f_k \Big )\leq 0   & \hspace{-0.4cm} $i \in \setN_g$ \label{ima_eq_g_c1} \\
- q_i + \mathbf{q}_i -   \mathbf{I}_{i,i} (e_i^2 + f_i^2)  -  \displaystyle{\sum_{k \in \delta(i)}} \Big [ \mathbf{I}_{i,k}(e_ie_k + f_if_k) +  \mathbf{R}_{i,k}(e_if_k - e_kf_i) \Big ]  -  \lambda(-A^3_i)  \displaystyle{\sum_{k\in \setN}} \Big  ( e_k^2 + f_k^2 - z^e_k - z^f_k \Big )\leq 0&  \hspace{-0.4cm} $i \in \setN_g$  \label{ima_eq_g_c2}  \\
- \mathbf{q}_i  +   \mathbf{I}_{i,i} (e_i^2 + f_i^2)  +  \displaystyle{\sum_{k \in \delta(i)}} \Big [  \mathbf{I}_{i,k}(e_ie_k + f_if_k) +   \mathbf{R}_{i,k}(e_if_k - e_kf_i) \Big ] -  \lambda(A^4_i) \displaystyle{\sum_{k\in \setN}}\Big  ( e_k^2 + f_k^2 - z^e_k - z^f_k \Big )\leq 0   & \hspace{-0.4cm}  $i \in \setN_d$  \label{ima_eq_d_c1}\\
\mathbf{q}_i -   \mathbf{I}_{i,i} (e_i^2 + f_i^2)  -  \displaystyle{\sum_{k \in \delta(i)}} \Big [  \mathbf{I}_{i,k}(e_ie_k + f_if_k) +   \mathbf{R}_{i,k}(e_if_k - e_kf_i) \Big ]-  \lambda(-A^4_i)  \displaystyle{\sum_{k\in \setN}} \Big  ( e_k^2 + f_k^2 - z^e_k - z^f_k \Big ) \leq 0  & \hspace{-0.4cm}  $i \in \setN_d$ \label{ima_eq_d_c2}
\end{numcases}

\end{scriptsize}

\medskip
By replacing Constraints~(\ref{real_eq_g})-(\ref{volt_magn}) by Constraints~(\ref{prod_e})-(\ref{volt_magn_c}),(\ref{real_eq_g_c1})-(\ref{ima_eq_d_c2}), we then obtain an equivalent problem to $(OPF)$ where the only non-convexity remains into the equalities (\ref{prod_e})-(\ref{prod_f}). Then, we use the McCormick's envelopes~\cite{Cor76}, to relax the latter equalities and we get a quadratic convex relaxation of $(OPF)$, that can then be used for the bounding step of a classical spatial \bb~algorithm. Performing such an algorithm is highly dependant on the quality of the bound at the root node. Moreover, we know that for the OPF problem the rank relaxation provides sharp bound. This is why in the rest of the paper, we focus on the computation of quadratic convex relaxation of $(OPF)$ whose value equals to the optimal value of the rank relaxation.

\bigskip

\section{Building a compact family of equivalent quadratic formulations}
\label{family}

\par We now present a compact family of quadratic reformulations to $(OPF)$. For simplicity, we start by rewriting the initial equality constraints~(\ref{real_eq_g})-(\ref{ima_eq_d}) by using the notation $y = (p,q) \in \mathbb{R}^{2\lvert \setN_g \rvert }$ and  $x = (e,f) \in \mathbb{R}^{2n}$ as follows:
 $$ \langle A_r , xx^T \rangle + a^T_ry   = b_r \quad r \in \setC  $$
where $\setC = (\setN_g, \setN_g, \setN_d,\setN_d)$, with $\lvert \setC\rvert =2n$, and $\forall \, r \in \setC$, $A_r \in \mathbf{S}_{2n}$ is the Hessian matrix of constraint $r$ (i.e. matrices $A^1_i$, $A^2_i$, $A^3_i$, and $A^4_i$), $a_r \in  \mathbb{R}^{2\lvert \setN_g \rvert}$ is the vector of linear coefficients of constraint $r$, and  $b \in \mathbb{R}^{2n}$, where coefficient $b_r$ is the the right-hand side of constraint $r$.

\medskip

\par Let $(\phi,\gamma) \in (\mathbb{R}^{2n},\mathbb{R}^{2n}) $ be two vector parameters, we build the following parameterized function:
\begin{align*}
  h_{\phi,\gamma}(x,y,z) = &    h(y)   +  \displaystyle{\sum_{r \in \setC}} \Big (  \phi_r( \langle A_r , xx^T \rangle + a^T_ry  - b_r)   \Big) +  \displaystyle{\sum_{i\in \setN}} \gamma_i \Big ( x_i^2 + x^2_{i+n} - z_i - z_{i+n} \Big )
\end{align*}

\noindent  where $h(y)$ is the initial objective function, and we recall that for all $i \in \setN$, $x_i = e_i$, $x_{i+n} = f_i$, $z_i = z^e_i$, and $z_{i+n} = z_i^f$. Observe that there exist parameters $(\phi,\gamma)$ such that  $h_{\phi,\gamma} $ is a convex function. Indeed, as mentioned above, function $h(y)$ is convex and separable. Now, the two additional terms are linear in $y$ and $z$, and pure quadratic in $x$. By taking $\forall \, r \in \setC$, $\bar{\phi}_r = 0$, and $\forall i \in \setN$, $\bar{\gamma}_i$ any non negative value,  the associated function $h_{\bar{\phi},\bar{\gamma}}(x,y,z)$ is obviously convex.

\medskip

By denoting $\lambda_r = \lambda_{min}(A_r)$ and  $\lambda_r' = \lambda_{min}(-A_r)$, we are now able to build a family of equivalent formulations to $(OPF)$:
\begin{numcases}{(OPF_{\phi,\gamma})} 
\min   h_{\phi,\gamma}(x,y,z)\nonumber \\
\mbox{s.t. } \nonumber\\
\quad  (\ref{real_bound}) (\ref{ima_bound})  (\ref{prod_e}) -(\ref{volt_magn_c}) \nonumber\\
\quad \langle A_r , xx^T \rangle + a^T_ry   -  \lambda_r\displaystyle{\sum_{i=1}^{2n } }(x_i^2 - z_i)   \leq b_r & $ r \in \setC$ \label{compact_ineg_1} \\
\quad \langle -A_r , xx^T \rangle - a^T_ry    -\lambda_r' \displaystyle{\sum_{i=1}^{2n } } (x_i^2 - z_i)  \leq - b_r & $ r \in \setC$ \label{compact_ineg_2}\\
\quad x=(e,f) \in \mathbb{R}^{2n}, \, \, y = (p,q) \in \mathbb{R}^{2\lvert \setN_g \rvert}, z=(z^e,z^f)  \in \mathbb{R}^{2n}  &  \label{setR}
\end{numcases}

where Constraints (\ref{compact_ineg_1})-(\ref{compact_ineg_2}) are Constraints (\ref{real_eq_g_c1})-(\ref{ima_eq_d_c2}) written in a compact form. It is easy to see that problem $(OPF_{\phi,\gamma})$ is equivalent to problem $(OPF)$, in the sense that any optimal solution from one is an optimal solution from the other. Indeed,  $h_{\phi,\gamma} (x,y,z) = h(y) $ when Constraints~(\ref{real_eq_g})-(\ref{ima_eq_d}), (\ref{prod_e}), and~(\ref{prod_f}) are satisfied.  Moreover, we already observed that we can choose parameters such that the objective function $h_{\phi,\gamma}$ is a convex function. Finally, the only constraints that remain non-convex are Constraints~(\ref{prod_e}) and~(\ref{prod_f}). To derive a quadratic convex relaxation of $(OPF_{\phi,\gamma})$, we relax them by use of the McCormick's upper envelopes (see~\cite{Cor76}), keeping the quadratic convex lower envelope. However, for this, we need upper and lower bounds on each variable $x_i$. Some trivial initial bounds can easily be deduced from Constraints (\ref{volt_magn}), i. e. $-\sqrt{\overline{\mathbf{v}}_i} \leq x_i \leq \sqrt{\overline{\mathbf{v}}_i}$, and $-\sqrt{\overline{\mathbf{v}}_i} \leq x_{i+n} \leq \sqrt{\overline{\mathbf{v}}_i}$.
\par By denoting with $(\ell,u) \in (\mathbb{R}^{2n}, \mathbb{R}^{2n})$ these bounds (i.e. $\ell = (- \sqrt{\overline{\mathbf{v}}}, - \sqrt{\overline{\mathbf{v}}})$ and $ u = ( \sqrt{\overline{\mathbf{v}}}, \sqrt{\overline{\mathbf{v}}})$), we replace Constraints~(\ref{prod_e}) and~(\ref{prod_f}) by the following set of convex inequalities:
\begin{numcases}{\setMC = (x,z)}
  z_i \leq   (u_i+\ell_i) x_i - u_i\ell_i & \nonumber\\
  z_i \geq   x^2_i  & \nonumber
 \end{numcases}

\medskip
We get  $(\overline{OPF}_{\phi,\gamma})$, a family of quadratic convex relaxations to $(OPF)$:
\begin{numcases}{(\overline{OPF}_{\phi,\gamma})} 
\min   h_{\phi,\gamma}(x,y,z)\nonumber \\
\mbox{s.t. } \nonumber\\
\quad  (\ref{real_bound}) (\ref{ima_bound})  (\ref{volt_magn_c}) \nonumber\\
\quad  (\ref{compact_ineg_1}) (\ref{compact_ineg_2}) (\ref{setR})  \nonumber\\
\quad (x,z) \in \setMC & \label{setMC}
\end{numcases}

\bigskip
\par Problem $(\overline{OPF}_{\phi,\gamma})$ is a compact quadratic convex relaxation to $(OPF)$, since we only add $\mathcal{O}(n)$ variables and constraints to the original formulation. 

\section{Computing a sharp quadratic convex relaxation}
\label{best_relax}
We are now interested in the best parameters $(\phi^*, \gamma^{*})$ that maximize the optimal value of $(\overline{OPF}_{\phi,\gamma})$ while making convex the parameterized function $h_{\phi,\gamma}$.   We prove that these best parameters can be deduced from the dual optimal solution of the rank relaxation of $(OPF)$.  For simplicity, we denote by $\gamma' = (\gamma,\gamma) \in \mathbb{R}^{2n}$, and  with this notation, we can rewrite function $h_{\phi,\gamma}$ as follows:
  \begin{align*}
  h_{\phi,\gamma'}(x,y,z) = &   \displaystyle{\sum_{i \in \setN_g}} \Big ( \mathbf{C}_{i,i} p^{2}_i  + \mathbf{c}_ip_i \Big )  + \langle \displaystyle{\sum_{r \in \setC}}\phi_r A_r  + \diag (\gamma'), xx^T \rangle +  \displaystyle{\sum_{r \in \setC}}  \phi_r( a^T_ry  - b_r)   -  \gamma^{'T}z
\end{align*}

\noindent where $\diag (\gamma')$ is the diagonal matrix of dimension $2n$, where the $i^{\textrm{th}}$-diagonal coefficient equals $\gamma'_i$. We formally pose the problem we aim to solve as follows:
\begin{align}{\OPT :}
 \max \Big \{   v(\overline{OPF}_{\phi,\gamma'})  \, : \, \displaystyle{\sum_{r \in \setC}}\phi_r A_r  + \diag (\gamma')  \succeq 0 \Big \} \nonumber
\end{align}

\noindent where  $v(\overline{OPF}_{\phi,\gamma'})$ is the optimal value of  problem $(\overline{OPF}_{\phi,\gamma'})$.
\medskip

We state in Theorem~\ref{theoopfquad} that the optimal value of  $\OPT$ equals the optimal value of the following the so-called rank relaxation of $(OPF)$:

\begin{numcases}{(SDP)} 
\min  h(Y,p) =  \displaystyle{\sum_{i \in \setN_g}} \Big ( \mathbf{C}_{i,i} Y_{i,i}  + \mathbf{c}_ip_i \Big ) \nonumber \\
\mbox{s.t. } \nonumber \\
\quad  \langle A_r , X \rangle + a^T_ry = b_r & $\forall r \in\setC $  \label{real_eq_g_sdp}\\
  \quad  X_{i,i} + X_{i+n,i+n} \leq \overline{\mathbf{v}}_i & $i \in \setN$ \label{volt_magn_sdp2} \\
 \quad - X_{i,i} - X_{i+n,i+n} \leq  -\underline{\mathbf{v}}_i& $i \in \setN$ \label{volt_magn_sdp1} \\
 \quad   \underline{\mathbf{p}}_i \leq p_i \leq  \overline{\mathbf{p}}_i & $i \in \setN_g$ \label{real_bound_sdp}\\
 \quad   \underline{\mathbf{q}}_i \leq q_i \leq  \overline{\mathbf{q}}_i & $i \in \setN_g$ \label{ima_bound_sdp}\\
\quad \left [ \begin{array}{ccc} 1 & y^T  \\
    y&  Y  \\ 
  \end{array} \right ] \succeq 0 \label{S_sdp1} \\
\quad X \succeq 0 \label{S_sdp2} \\
\quad  y \in \mathbb{R}^{2\lvert \setN_g \rvert},    (Y,X)  \in (\mathbf{S}_{2\lvert \setN_g \rvert} ,\mathbf{S}_{2n})&  \nonumber
\end{numcases}

\noindent where Constraints~(\ref{real_eq_g_sdp}) are the compact form of the following set of constraints:

\begin{small}
\begin{numcases}{} 
 - p_i + \mathbf{R}_{i,i} (X_{i,i} + X_{i+n,i+n})  -  \displaystyle{\sum_{k \in \delta(i)}} \Big [ \mathbf{R}_{i,k}(X_{i,k} + X_{i+n,k+n}) -   \mathbf{I}_{i,k}(X_{i,k+n} - X_{k,i+n}) \Big ] = - \mathbf{p}_i &$i \in \setN_g$ \nonumber \\
 \mathbf{R}_{i,i} (X_{i,i} + X_{i+n,i+n})  -  \displaystyle{\sum_{k \in \delta(i)}} \Big [ \mathbf{R}_{i,k}(X_{i,k} + X_{i+n,k+n}) -   \mathbf{I}_{i,k}(X_{i,k+n} - X_{k,i+n}) \Big ]= - \mathbf{p}_i&$i \in \setN_d$ \nonumber \\
 - q_i +  \mathbf{I}_{i,i} (X_{i,i} + X_{i+n,i+n})  -  \displaystyle{\sum_{k \in \delta(i)}} \Big [ - \mathbf{I}_{i,k}(X_{i,k} + X_{i+n,k+n}) -   \mathbf{R}_{i,k}(X_{i,k+n} - X_{k,i+n}) \Big ]= - \mathbf{q}_i&$i \in \setN_g$ \nonumber \\
 \mathbf{I}_{i,i} (X_{i,i} + X_{i+n,i+n})  -  \displaystyle{\sum_{k \in \delta(i)}} \Big [ - \mathbf{I}_{i,k}(X_{i,k} + X_{i+n,k+n}) -   \mathbf{R}_{i,k}(X_{i,k+n} - X_{k,i+n}) \Big ]= - \mathbf{q}_i&$i \in \setN_d$ \nonumber 
\end{numcases}
\end{small}

\medskip

\begin{theorem}
  \label{theoopfquad}
 The optimal value of  $\OPT $ equals the optimal value of $(SDP)$.
\end{theorem}

\begin{preuve} 

  \noindent $\diamond$ To prove that $v\OPT \leq v(SDP)$, we show that  for any feasible solution $(\bar{\phi},\bar{\gamma}')$ to $(\overline{OPF}_{\bar{\phi} ,\bar{\gamma}'})$, we have $v(\overline{OPF}_{\bar{\phi},\bar{\gamma}'}) \leq v(SDP)$, which in turn implies that $v\OPT \leq v(SDP)$ since the right hand side is constant. Let $(\bar{y},\bar{Y},\bar{X})$ be a feasible solution of $(SDP)$, and build the solution $(x =0,y = \bar{y},z = \Diag (\bar{X}))$, where $\Diag (X)$ is the vector composed of the diagonal terms of matrix $X$. We show that: i)  $(x,y,z)$ is feasible for $(\overline{OPF}_{\bar{\phi} ,\bar{\gamma}'})$, and i,i) its objective value is less or equal than $v(SDP)$. Since $(\overline{OPF}_{\bar{\phi} ,\bar{\gamma}'})$ is a minimization problem, $(\overline{OPF}_{\bar{\phi} ,\bar{\gamma}'}) \leq v(SDP)$ follows.
  
\medskip
  
\begin{enumerate}[i)]
\item   We show that $(0,\bar{y},\Diag (\bar{X}))$ is  feasible  to $(\overline{OPF}_{\bar{\phi} ,\bar{\gamma}'})$.  Obviously,  Constraints~(\ref{real_bound})-(\ref{ima_bound}), and~(\ref{volt_magn_c}) are satisfied. We prove now prove that Constraints~(\ref{compact_ineg_1}), (\ref{compact_ineg_2}), and (\ref{setMC})  are satisfied. 
  \begin{enumerate}
  \item \textit{Constraints~(\ref{compact_ineg_1})}. We start by observing that, since $x=0$, $\lambda_r \leq 0$ and $\bar{X}_{i,i} \geq 0 $, we have for all $r \in~\setC$, $-\lambda_r \displaystyle{\sum_{i=1}^{2n } } (x_i^2 - z_i) = -\lambda_r \displaystyle{\sum_{i=1}^{2n } } (x_i^2 - \bar{X}_{i,i}) = \lambda_r \displaystyle{\sum_{i=1}^{2n } }( \bar{X}_{i,i})  \leq 0 $. Constraints~(\ref{compact_ineg_1}) can be rewritten as:
\begin{align*}
     b_r  - a^T_r\bar{y}  - \langle A_r , \bar{x}\bar{x}^T \rangle +\lambda_r \displaystyle{\sum_{i=1}^{2n } } (\bar{x}_i^2 - \bar{X}_{i,i})   =  b_r  - a^T_r\bar{y} -  \lambda_r\displaystyle{\sum_{i=1}^{2n } } \bar{X}_{i,i} \geq 0\\
    \shortintertext{since  $- \lambda_r\displaystyle{\sum_{i=1}^{2n } }\bar{X}_{i,i} \geq 0$ and by Constraints~(\ref{real_eq_g_sdp}) and (\ref{S_sdp2}), $b_r  - a^T_r\bar{y}  = \langle A_r , \bar{X} \rangle \geq 0 $}
    \end{align*}
\item \textit{Constraints~(\ref{compact_ineg_2})}. Similarly, Constraints~(\ref{compact_ineg_2}) can be rewritten as:  
     \begin{align*}
       -b_r  + a^T_r\bar{y}  + \langle A_r , \bar{x}\bar{x}^T \rangle +  \lambda_r' \displaystyle{\sum_{i=1}^{2n } }(\bar{x}_i^2 - \bar{X}_{i,i})   =  -b_r  + a^T_r\bar{y}  -  \lambda_r'\displaystyle{\sum_{i=1}^{2n } } \bar{X}_{i,i} \geq 0\\
       \shortintertext{since  $- \lambda_r'\displaystyle{\sum_{i=1}^{2n } }\bar{X}_{i,i} \geq 0$ and by Constraints~(\ref{real_eq_g_sdp}) and (\ref{S_sdp2}), $-b_r  + a^T_r\bar{y}  = \langle -A_r , \bar{X} \rangle \geq 0 $}
     \end{align*}
   \item \textit{Constraints~(\ref{setMC})}. Since $x =0$, the set $\setMC$ becomes:
     \begin{numcases}{}
       \bar{X}_{i,i} \leq   \overline{\mathbf{v}}_i   & \nonumber\\
       \bar{X}_{i,i} \geq 0 & \nonumber
     \end{numcases}
Since $X \succeq 0$, we have  $\bar{X}_{i,i} \geq 0$. Moreover,  combining $\bar{X}_{i+n,i+n} \geq 0$  with Constraints~(\ref{volt_magn_sdp2}) we get the first inequality.
     
 \end{enumerate}

\item  Let us now compare the objective values: we prove that $\Delta = h_{\phi,\gamma'}(y,0,z) - h(Y,p) \leq 0$.
     \begin{align*}
      \Delta= &  \displaystyle{\sum_{i \in \setN_g}}  \mathbf{C}_{i,i} \bar{y}^{2}_i    +  \displaystyle{\sum_{r \in \setC}}  \bar{\phi}_r( a^T_r\bar{y}  - b_r)    +  \displaystyle{\sum_{i \in \setN}} \bar{\gamma}_i  ( - \bar{z}^1_{i} - \bar{z}^2_{i}  )  -  \displaystyle{\sum_{i \in \setN_g}}\mathbf{C}_{i,i} \bar{Y}_{i,i}   \\
      \Delta= &  \displaystyle{\sum_{i \in \setN_g}} \mathbf{C}_{i,i} (\bar{y}^{2}_i - \bar{Y}_{i,i} )+  \displaystyle{\sum_{r \in \setC}} \bar{\phi}_r( a^T_r\bar{y}  - b_r)  +  \displaystyle{\sum_{i \in \setN}} \bar{\gamma}_i ( - \bar{X}_{i,i} - \bar{X}_{i+n,i+n}  )\\
      \shortintertext{By constraints~(\ref{real_eq_g_sdp}), an by definition of $\gamma'$, we have:}
      \Delta= &  \displaystyle{\sum_{i \in \setN_g}} \mathbf{C}_{i,i} (\bar{y}^{2}_i - \bar{Y}_{i,i} ) -   \displaystyle{\sum_{r \in \setC}}\bar{\phi}_r\langle A_r, \bar{X}\rangle - \langle  \diag (\bar{\gamma}'), \bar{X} \rangle \\
      \Delta= &  \displaystyle{\sum_{i \in \setN_g}} \mathbf{C}_{i,i} (\bar{y}^{2}_i - \bar{Y}_{i,i} )-   \displaystyle{\sum_{r \in \setC}}\bar{\phi}_r\langle A_r  + \diag (\bar{\gamma}') , \bar{X}\rangle \\
       \Delta \leq  &0\\
        \shortintertext{Indeed, the first term is non positive since $\mathbf{C}_{i,i}\geq 0$, and by Constraint~(\ref{S_sdp1}), we have that $\forall i, \, \,  \bar{y}^{2}_i - \bar{Y}_{i,i} \leq 0$. This is also the case for the last term, by   Constraint~(\ref{S_sdp2}), and since $ \displaystyle{\sum_{r \in \setC}}\bar{\phi}_r A_r + \diag (\bar{\gamma}')  \succeq 0$ by feasibility of  $\OPT$.}
  \end{align*}
\end{enumerate}
\noindent  $\diamond$ Let us secondly prove that $v\OPT \geq v(SDP)$ or equivalently  $v\OPT \geq v(D)$ where $(D)$ is the dual of $(SDP)$:
\begin{numcases}{(D)}
  \max \displaystyle{\sum_{i \in \setN_g}}( \underline{\theta}_i^p \underline{\mathbf{p}}_i + \underline{\theta}_i^q \underline{\mathbf{q}}_i   -\overline{\theta}_i^p\overline{\mathbf{p}}_i - \overline{\theta}_i^q\overline{\mathbf{q}}_i  ) +\displaystyle{\sum_{i \in \setN}}(  \underline{\gamma}_i \underline{\mathbf{v}}_i- \overline{\gamma}_i \,\overline{\mathbf{v}}_i )  - \displaystyle{\sum_{r \in \setC}}\phi_r b_r   - \rho \nonumber\\
 \mbox{s.t.} \nonumber \\
   W=\left [ \begin{array}{ccc}
 \rho &  \frac{1}{2} (c +  \overline{\theta} -  \underline{\theta} + \displaystyle{\sum_{r \in \setC}}\phi_r a_{r})^T &  \mathbf{0}^T_{2n} \\
 \frac{1}{2} (c +  \overline{\theta} -  \underline{\theta} + \displaystyle{\sum_{r \in \setC}}\phi_r a_{r})&    \mathbf{C} &  \mathbf{0}_{2 \lvert \setN_g \rvert, 2n}  \\
  \mathbf{0}_{2n} & \mathbf{0}_{2n,2 \lvert \setN_g \rvert}   &  \displaystyle{\sum_{r \in \setC}}\phi_r A_r + \diag (\gamma')
   \end{array} \right ]\succeq 0  \label{dsdp_1}  \\
  \quad  \overline{\theta} =  (\overline{\theta}^p,\overline{\theta}^q) ,  \underline{\theta} =  (\underline{\theta}^p,\underline{\theta}^q), \gamma' =  ((\overline{\gamma} - \underline{\gamma}), (\overline{\gamma} - \underline{\gamma})) \label{dsdp_2}  \\
  \quad C = \left [ \begin{array}{cc} C& \mathbf{0}_{ \lvert \setN_g \rvert,\lvert \setN_g \rvert}\\
   \mathbf{0}_{ \lvert \setN_g \rvert,\lvert \setN_g \rvert} &  \mathbf{0}_{ \lvert \setN_g \rvert,\lvert \setN_g \rvert}\\
    \end{array} \right ] \in \mathbf{S}_{2 \lvert \setN_g \rvert} \label{dsdp_3}  \\
  \quad c =(\mathbf{c}, \mathbf{0}_{ \lvert \setN_g \rvert + 2n}) \label{dsdp_4}  \\
 \quad (\phi,\overline{\gamma}, \underline{\gamma},  \overline{\theta}^p, \underline{\theta}^p,\overline{\theta}^q, \underline{\theta}^q,\rho ) \in (\mathbb{R}^{|\setC|},\mathbb{R}^{n},\mathbb{R}^{n}, \mathbb{R}^{|\setN_g|}, \mathbb{R}^{|\setN_g|},\mathbb{R}^{|\setN_g|},\mathbb{R}^{|\setN_g|},\mathbb{R})& \nonumber
\end{numcases}

\noindent where  $\phi, \overline{\gamma}$, $\underline{\gamma}$, $(\overline{\theta}^p, \underline{\theta}^p)$, $(\overline{\theta}^q, \underline{\theta}^q)$, $\rho$  are the dual variables associated to Constraints~(\ref{real_eq_g_sdp})-(\ref{S_sdp1}) respectively. We denote by $\mathbf{0}_{n} $ ( $\mathbf{0}_{n,n} $ respectively ) the $n$-dimensional ($n\times n$-dimensional resp.) vector (matrix resp. ) where each coefficient equals $0$.

\par Let $ (\phi^*,\overline{\gamma}^*, \underline{\gamma}^*,  \overline{\theta}^{p*}, \underline{\theta}^{p*},\overline{\theta}^{q*}, \underline{\theta}^{q*},\rho^* ) $ be an optimal solution to $(D)$, we build the following solution $(\phi^* = \phi^*, \gamma^{'*} = ( (\overline{\gamma}^* - \underline{\gamma}^* ) ,(\overline{\gamma}^* - \underline{\gamma}^* ))$ that is feasible for $\OPT$, i.e. $\displaystyle{\sum_{r \in \setC}}\phi^*_r A_r + \diag (\gamma^{'*}) \succeq 0$, by Constraint~(\ref{dsdp_1}). The objective value of this solution is equal to  $v(\overline{OPF}_{\phi^*,\gamma^{'*}})$. Hence to prove that $v(\overline{OPF}_{\phi^*,\gamma^{'*}}) \geq v(D)$,  we prove that for any feasible solution $(x,y,z)$ to $(\overline{OPF}_{\phi^*,\gamma^{'*}}) $, the associated objective value is not smaller than the optimal value of $(D)$. Denote by  $\Delta$ the difference between the objective values, and we below prove that $\Delta \geq 0$.
\begin{align*}
  \Delta & = \displaystyle{\sum_{i \in \setN_g}} \Big ( \mathbf{C}_{i,i} p^{2}_i  + \mathbf{c}_ip_i \Big )  + \langle \displaystyle{\sum_{r \in \setC}}\phi^*_r A_r  + \diag (\gamma^{'*}), xx^T \rangle +  \displaystyle{\sum_{r \in \setC}}  \phi^*_r( a^T_ry  - b_r)  -  \gamma^{'*T}z  \\
  & \quad \quad - \displaystyle{\sum_{i \in \setN_g}}( \underline{\theta}_i^{p*} \underline{\mathbf{p}}_i  + \underline{\theta}_i^{q*} \underline{\mathbf{q}}_i   -\overline{\theta}_i^{p*}\overline{\mathbf{p}}_i - \overline{\theta}_i^{q*}\overline{\mathbf{q}}_i ) -\displaystyle{\sum_{i \in \setN}}(\underline{\gamma}_i^* \underline{\mathbf{v}}_i - \overline{\gamma}_i^* \,\overline{\mathbf{v}}_i  )  + \displaystyle{\sum_{r \in \setC}}\phi^*_r b_r   + \rho^*  \\
  \Delta & =\langle C, yy^T \rangle + \langle \displaystyle{\sum_{r \in \setC}}\phi^*_r A_r  + \diag (\gamma^{'*}), xx^T \rangle  - \displaystyle{\sum_{i \in \setN_g}}( \underline{\theta}_i^{p*} \underline{\mathbf{p}}_i  + \underline{\theta}_i^{q*} \underline{\mathbf{q}}_i   -\overline{\theta}_i^{p*}\overline{\mathbf{p}}_i - \overline{\theta}_i^{q*}\overline{\mathbf{q}}_i ) \\
  & \quad \quad +  \displaystyle{\sum_{r \in \setC}}  \phi^*_r a^T_ry  +  \displaystyle{\sum_{i \in \setN_g}}\mathbf{c}_ip_i  -\displaystyle{\sum_{i \in \setN}} (\overline{\gamma}_i^*  - \underline{\gamma}_ i^* ) ( z_i^e  + z_i^f) -\displaystyle{\sum_{i \in \setN}}(\underline{\gamma}_i^* \underline{\mathbf{v}}_i - \overline{\gamma}_i^* \,\overline{\mathbf{v}}_i  )  + \rho^*  \\
  \intertext{By Constraints~(\ref{real_bound_sdp}) and (\ref{ima_bound_sdp}) and since all coefficients $\overline{\theta}_i^{p*}, \, \underline{\theta}_i^{p*}, \, \overline{\theta}_i^{q*}, \, \underline{\theta}_i^{q*}$ are non-negative, we get:}
  \Delta & \geq\langle C, yy^T \rangle + \langle \displaystyle{\sum_{r \in \setC}}\phi^*_r A_r  + \diag (\gamma^{'*}), xx^T \rangle  - \displaystyle{\sum_{i \in \setN_g}}( \underline{\theta}_i^{p*} p_i + \underline{\theta}_i^{q*} q_i   -\overline{\theta}_i^{p*}p_i - \overline{\theta}_i^{q*}q_i ) \\
  & \quad \quad +  \displaystyle{\sum_{r \in \setC}}  \phi^*_r a^T_ry  + c^Ty + \displaystyle{\sum_{i \in \setN}} \big ( \overline{\gamma}_i^*(\overline{\mathbf{v}}_i   - z_i^e  - z_i^f) + \underline{\gamma}_ i^* ( z_i^e  + z_i^f - \underline{\mathbf{v}}_i ) \big ) + \rho^*  \\
\intertext{By Constraints~(\ref{volt_magn_c})  and since all coefficients $\overline{\gamma}_i^*, \, \underline{\gamma}_i^*$ are non-negative, we get:}
\Delta & \geq\langle C, yy^T \rangle + \langle \displaystyle{\sum_{r \in \setC}}\phi^*_r A_r  + \diag (\gamma^{'*}), xx^T \rangle  +  (\overline{\theta}- \underline{\theta})^T y +  \displaystyle{\sum_{r \in \setC}}  \phi^*_r a^T_ry  + c^Ty + \rho^*  \\
 \Delta & \geq\langle C, yy^T \rangle + \langle \displaystyle{\sum_{r \in \setC}}\phi^*_r A_r  + \diag (\gamma^{'*}), xx^T \rangle  +  (\overline{\theta}- \underline{\theta} +  \displaystyle{\sum_{r \in \setC}}  \phi^*_r a_r + c)^Ty + \rho^*  
\end{align*}
Finally, by Constraint~(\ref{dsdp_1}), we have  \begin{small}$ \left ( \begin{array}{ll}    1 \\   y  \\x    \end{array}\right )^T W   \left ( \begin{array}{ll}    1 \\   y   \\ x  \end{array}\right ) \geq 0$ \end{small}, which prove that $\Delta \geq 0$.

\end{preuve}

From the proof of Theorem~\ref{theoopfquad},  we can characterize the best parameters $\phi$ and $\gamma$

\begin{corollary}
  The best parameters $(\phi^*,\gamma^{'*} = \Big ( (\overline{\gamma}^* - \underline{\gamma}^*),(\overline{\gamma}^* - \underline{\gamma}^*) \Big )$ can be computed by solving $(SDP)$. In particular, $\phi^*$ is the vector of optimal dual variables associated to Constraints~(\ref{real_eq_g_sdp}), and  $\overline{\gamma}^*, \underline{\gamma}^*$, the vectors of optimal dual variables associated to Constraints~(\ref{volt_magn_sdp2}) and~(\ref{volt_magn_sdp1}), respectively.
\end{corollary}

Hence, for computing the best convex relaxation of the OPF problem, we need to solve once problem $(SDP)$.  Then, to solve $(OPF)$ to global optimality, we perform a spatial \bb~algorithm based on the bound given by problem $(\overline{OPF}_{\phi^*,\lambda^*})$. As proven in Theorem~\ref{theoopfquad}, the bound at the root node of the tree equals the optimal value of $(SDP)$. Another advantage of our compact relaxation, in comparison to the complete linearization of the constraints of method \rcopf, is that we only have to enforce $2n$ equalities to prove global optimality in our \bb. We sum up our global optimization algorithm \texttt{Compact OPF} (\opfdiag) in Algorithm~\ref{algo_opfdiag}.

\begin{algorithm}
\caption{Solution algorithm \opfdiag~for exact solution of $(OPF)$}

\label{algo_opfdiag}

\begin{algorithmic}
\STATE  \begin{enumerate}[{\textbf{step}} 1\textbf{:}]
\item Solve $(SDP)$.
\vspace{0.2cm}
\item Deduce $\phi^*$ and $\gamma^{*'}$.
\vspace{0.2cm}
\item Solve $(OPF_{\phi^*,\gamma^{*'}})$ with a \sbb~based on the quadratic convex relaxation $(\overline{OPF}_{\phi^*,\gamma^{*'}})$ .
\end{enumerate}
\end{algorithmic}

\end{algorithm}



\section{Some experimental results}
\label{result}

\par In this section, we illustrate on some experiments the behavior of the algorithm \opfdiag~for exact solution of instances of the \textit{PG-lib} library~\cite{PGLIB}. We compare \opfdiag~with the non-linear solver  \baron~\cite{baron}, and with the initial quadratic convex reformulation approach \rcopf~\cite{Godard19, EGLMR19}. Our experiments were carried out on a server with $2$ CPU Intel Xeon each of them having $12$ cores and $2$ threads of $2.5$ GHz and $4*16$ GB of RAM using a Linux operating system. We set the time limit to 3 hours for all methods. For the solver \baron, we use the multi-threading version of \texttt{Cplex 12.9}~\cite{cplex129}  with up to 64 threads. For methods \opfdiag~and \rcopf, we used the semidefinite solver \mosek~\cite{mosek} for solving semidefinite programs. At each node of the spatial branch-and-bound, we used the solver  \mosek~for solving the QCQP of method \opfdiag, and the solver \texttt{Cplex 12.9} for solving QP of method \rcopf. For computing feasible local solutions, we use the local solver \ipopt~\cite{ipopt}. 

\bigskip
\par For our experiences, we considered medium-sized data of power networks having 3 to 300 buses, and we took the  formulation of the OPF described by Constraints~(\ref{real_eq_g})-(\ref{ima_bound}). We report in Table~1 the characteristics of each instance: its \textit{Name}, and the number of  \textit{Buses}, \textit{Generators}, and \textit{Lines} of the considered  power network. We also indicate in the column \textit{Opt} the best solution found by \opfdiag~and \rcopf, within 3 hours of computing time. The column $\lvert(y,x)\rvert$ specifies the dimension of the variable vector $(y,x)$ in \OPF.

\begin{table}
\centering
\begin{small}
  \begin{tabular}{|l|c|c|c|c|c|} \hline
\textit{Name}&\textit{Buses}&\textit{Generators}&\textit{Lines}&\textit{Opt}&$\lvert(x,y)\rvert$\\\hline
\texttt{caseWB2}&2&1&1&9.0567&6\\\hline
\texttt{caseWB3}&3&1&2&417.2453&8\\\hline
\texttt{pglib\_opf\_case3\_lmbd}&3&3&3&5 694.5249&12\\\hline
\texttt{caseWB5}&5&2&6&13.7797&14\\\hline
\texttt{pglib\_opf\_case5\_pjm}&5&5&5&14 997.0431&20\\\hline
\texttt{case6ww}&6&3&11&3 126.3145&18\\\hline
\texttt{pglib\_opf\_case14\_ieee}&14&5&20&2 178.0893&38\\\hline
\texttt{pglib\_opf\_case24\_ieee\_rts}&24&33&38&63 344.6382&114\\\hline
\texttt{pglib\_opf\_case30\_as}&30&6&41&801.5451&72\\\hline
\texttt{pglib\_opf\_case30\_ieee}&30&6&41& 6 592.9534&72\\\hline
\texttt{pglib\_opf\_case39\_epri}&39&10&46&133 801.7063&98\\\hline
\texttt{pglib\_opf\_case57\_ieee}&57&7&80&37 589.3248 &128\\\hline
\texttt{pglib\_opf\_case73\_ieee\_rts}&73&99&120&189 741.3755  &344\\\hline
\texttt{pglib\_opf\_case89\_pegase}&89&12&210&106 696.9325 &202\\\hline
\texttt{pglib\_opf\_case118\_ieee}&118&54&186&96 881.5257&344\\\hline
\texttt{pglib\_opf\_case162\_ieee\_dtc}&162&12&284&84 785.2377&348\\\hline
\texttt{pglib\_opf\_case179\_goc}&179&29&263&750 158.5809&416\\\hline
\texttt{pglib\_opf\_case200\_activ}&200&38&245& 27 557.5673&476\\\hline
\texttt{pglib\_opf\_case240\_pserc}&240&143&448& - &766\\\hline
\texttt{pglib\_opf\_case300\_ieee}&300&69&411&- &738\\\hline
\end{tabular}
  \end{small}
\label{tab:carac}
\caption {Characteristics of the considered instances of \textit{PG-lib} library.}

\end{table}

\bigskip

\par In Figure~1, we present the performance profile of the CPU times for methods \opfdiag, \rcopf, and the solver \texttt{Baron 19.3.24}. We observe that  \opfdiag~and \rcopf~significantly outperform the solver \baron. In fact \baron~solves to optimality only 6 instances out of the 20 considered, the largest of which is \texttt{pglib\_opf\_case14\_ieee}. The two other approaches are more efficient, since they solve within 3 hours of CPU time, 12 instances for \rcopf~and 14 for \opfdiag. Moreover,  this profile shows that \opfdiag~is faster than \rcopf~for these instances.

\begin{figure}
\begin{center}
 \includegraphics[width=12cm]{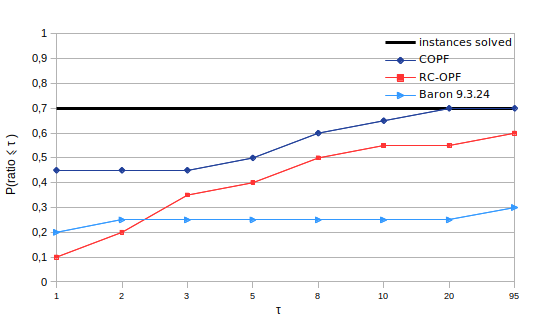}
\label{perfprof}
 \caption{Performance profile of the total time for networks with $2$ to $300$ buses (time limit 3 hours). }
\end{center}
\end{figure}
\bigskip

\begin{table}
\centering
\begin{small}
  \begin{tabular}{|l|c||c|c|c|c|c||c|c|c|c|c|} \hline
 \multicolumn{2}{|c||}{Instance} & \multicolumn{5}{c||}{\opfdiag}& \multicolumn{5}{c|}{\rcopf}\\\hline
 \textit{Name}&\textit{Gap}&$\lvert z \rvert $&\textit{UB}&\textit{SDP}&\textit{CPU}&\textit{Nodes} &$\lvert z \rvert $&\textit{UB}&\textit{SDP}&\textit{CPU}&\textit{Nodes} \\\hline
\texttt{caseWB2}&1.947&4&0&1&2&25&8&0&1&-&44265\\\hline
\texttt{caseWB3}&0.000&6&0&1&1&1&14&0&1&1&0\\\hline
\texttt{pglib\_opf\_case3\_lmbd}&0.000&6&0&1&1&1&18&0&1&1&0\\\hline
\texttt{caseWB5}&28.093&10&0&1&357&8 267&34&0&1&-&292231\\\hline
\texttt{pglib\_opf\_case5\_pjm}&0.000&10&0&1&1&0&34&0&1&1&0\\\hline
\texttt{case6ww}&0.000&12&3&1&4&0&56&0&1&1&0\\\hline
\texttt{pglib\_opf\_case14\_ieee}&0.000&28&0&1&1&0&108&0&1&2&0\\\hline
\texttt{pglib\_opf\_case24\_ieee\_rts}&0.000&48&0&1&1&0&184&42&1&43&0\\\hline
\texttt{pglib\_opf\_case30\_as}&0.000&60&8&0&8&0&224&0&1&1&0\\\hline
\texttt{pglib\_opf\_case30\_ieee}&0.000&60&0&1&1&0&224&5&1&6&0\\\hline
\texttt{pglib\_opf\_case39\_epri}&0.000&78&1&1&2&0&262&5&1&6&0\\\hline
\texttt{pglib\_opf\_case57\_ieee}&0.003&114&167&1& - &2 659&426&31&1& -&2193\\\hline
\texttt{pglib\_opf\_case73\_ieee\_rts}&0.000&146&67&2&70&0&578&191&3&196&0\\\hline
\texttt{pglib\_opf\_case89\_pegase}&0.000&178&31&3&34&0&1002&72&3&76&0\\\hline
\texttt{pglib\_opf\_case118\_ieee}&0.005&236&933&2& - &973&952&2864&5& - &1023\\\hline
\texttt{pglib\_opf\_case162\_ieee\_dtc}&1.890&324&342&4&-&353&1444&518&7& -&1033\\\hline
\texttt{pglib\_opf\_case179\_goc}&0.034&358&1517&5& - &343&1246&3045&10& -&1127\\\hline
\texttt{pglib\_opf\_case200\_activ}&0.000&400&420&5&426&0&1380&2156&10&2172&0\\\hline
\texttt{pglib\_opf\_case240\_pserc}&-&480&3624&27&-&219&1872&3328&27&-&619\\\hline
\texttt{pglib\_opf\_case300\_ieee}&-&600&3295&42&-&29&2236&3124&41&-&1171\\\hline
\end{tabular}
  \end{small}
\label{tab:results}
\caption {Initial gap, Sizes, CPU times and Nodes for methods \opfdiag~and \rcopf.}

\end{table}
\par Finally, we present in Table~2, a detailed comparison between methods \opfdiag~and \rcopf.  Column \textit{Gap}: $ = \displaystyle \left|\frac {Opt - Cont }{Opt} \right| * 100$, is the initial gap at the root node of the \bb, where \textit{Cont} is the optimal value of the rank relaxation,  \textit{Opt} is defined as in Table~1, and $-$ indicates that no feasible solution has been found by the algorithm. Column $\lvert z \rvert $ specifies  the number auxiliary variables in the relaxation, that is also the number of non-convex equalities to force during the spatial \bb. Columns \textit{UB} and \textit{SDP} report the CPU times in seconds for finding an initial upper bound and for solving the rank relaxation. Column \textit{CPU} is the total CPU time, where $-$ means that the instance is unsolved within the time limit of 3 hours. Finally, Column \textit{Nodes} is the number of nodes visited by the \bb. A first observation concerns the reformulation time which is significantly shorter than the global resolution time (always less than 42 seconds), while the rank relaxation is solved directly with a standard semidefinite solver. The zero gap at the root of the \bb~for 12 instances out of the 20 considered confirms the strength of the rank relaxation for the OPF problem. On the other hand, these experiments clearly show that the exact resolution of the instances where the gap is non-zero remains very difficult. Indeed, the method \rcopf~does not solve instances having only 2 or 5 buses in the considered power network. In fact, during its \bb, \rcopf~does not increase the lower bound, despite a large number of explored nodes. This is not the case for \opfdiag, which, with a much smaller number of nodes, slightly increases the lower bound over the course of the \bb, even for the largest instances. This is because the number of auxilliary variables and relaxed equalities $z_i = x^2_i$  is  strongly reduced in \opfdiag. Let us finally note that determining a feasible solution is also difficult in practice, at least with a generic interior point algorithm. Moreover, the CPU time for this step is not constant since  for example it goes from 400 to 1380 seconds for the instance \texttt{pglib\_opf\_case200\_activ}.

\section{Conclusion}

\label{conc}
We consider the OPF problem that determines the power production at each bus  of an electric network minimizing a production cost. In this paper, we introduce a global optimisation algorithm that is based on a new quadratically constrained quadratic relaxation. This relaxation is compact in the sense that it has only $\mathcal{O}(n)$ auxiliary variables and constraints, where $n$ is the number of buses of the network. We moreover prove that our quadratic relaxation has the same optimal value as the rank relaxation. Finally, to solve $(OPF)$ to global optimality, we perform a spatial \bb~algorithm based on our new quadratic convex relaxation. Another advantage of our approach is that to prove global optimality, we have a reduced number of non-convex equalities to force in spatial the \bb. We report computational results on instances of the literature. These results show that this new approach is more efficient than the standard solver \baron. Moreover, for the most difficult instances, and with a basic implementation, it is able to improve the lower bound unlike other methods in the literature. A future work consists in using OBBT techniques to further improve the behaviour of \opfdiag~on the most difficult instances.

\bigskip
\bibliographystyle{plain}
\bibliography{mybib}

\begin{thebibliography}{10}

\bibitem{mosek}
MOSEK ApS.
\newblock {\em The MOSEK optimization toolbox for MATLAB manual. Version 9.2.},
  2019.

\bibitem{BW09}
X.~Bai and H.~Wei.
\newblock Semi-definite programming-based method for security-constrained unit
  commitment with operational and optimal power flow constraints.
\newblock {\em IET Generation, Transmission \& Distribution}, 3:182--197(15),
  February 2009.

\bibitem{BWFW08}
X.~Bai, H.~Wei, K.~Fujisawa, and Y.~Wang.
\newblock Semidefinite programming for optimal power flow problems.
\newblock {\em International Journal of Electrical Power \& Energy Systems},
  30(6):383 -- 392, 2008.

\bibitem{CHV15b}
C.~{Coffrin}, H.~L. {Hijazi}, and P.~{Van Hentenryck}.
\newblock The qc relaxation: A theoretical and computational study on optimal
  power flow.
\newblock {\em IEEE Transactions on Power Systems}, 31(4):3008--3018, 2016.

\bibitem{CHV17}
C.~Coffrin, H.L. Hijazi, and P.~Van Hentenryck.
\newblock Convex quadratic relaxations for mixed-integer nonlinear programs in
  power systems.
\newblock {\em Mathematical Programming Computation}, 9(3):321--367, 2017.

\bibitem{CofVan14}
Carleton Coffrin and Pascal~Van Hentenryck.
\newblock A linear-programming approximation of ac power flows.
\newblock {\em INFORMS J. Comput.}, 26:718--734, 2014.

\bibitem{EGLMR19}
S.~Elloumi, H.~Godard, A.~Lambert, J.~Maegth, and M.~Ruiz.
\newblock Global optimality of optimal power flow using quadratic convex
  optimization.
\newblock {\em 6th {I}nternational {C}onference on {C}ontrol, {D}ecision and
  {I}nformation {T}echnologies, CODIT}, pages 1--6, 2019.

\bibitem{EllLam19}
S.~Elloumi and A.~Lambert.
\newblock Global solution of non-convex quadratically constrained quadratic
  programs.
\newblock {\em Optimization Methods and Software}, 34(1):98--114, 2019.

\bibitem{GLTL15}
L.~{Gan}, N.~{Li}, U.~{Topcu}, and S.~H. {Low}.
\newblock Exact convex relaxation of optimal power flow in radial networks.
\newblock {\em IEEE Transactions on Automatic Control}, 60(1):72--87, 2015.

\bibitem{Godard19}
H.~Godard.
\newblock {\em R\'{e}solution exacte du probl\`{e}me de l'optimisation des flux
  de puissance}.
\newblock Th\`{e}se de doctorat en informatique, Conservatoire National des
  Arts et M\'{e}tiers, Paris, 2019.

\bibitem{HamRub70}
P.L. Hammer and A.A. Rubin.
\newblock Some remarks on quadratic programming with 0-1 variables.
\newblock {\em Revue Fran\c{c}aise d'Informatique et de Recherche
  Op\'erationnelle}, 4:67--79, 1970.

\bibitem{cplex129}
IBM.
\newblock {IBM CPLEX 12.9 R}eference {M}anual.
\newblock
  "\url{https://www.ibm.com/docs/en/icos/12.9.0?topic=cplex-callable-library-c-api-reference-manual}",
  2020.

\bibitem{Jab06}
R.~A. {Jabr}.
\newblock Radial distribution load flow using conic programming.
\newblock {\em IEEE Transactions on Power Systems}, 21(3):1458--1459, 2006.

\bibitem{Jab12}
R.~A. {Jabr}.
\newblock Exploiting sparsity in sdp relaxations of the opf problem.
\newblock {\em IEEE Transactions on Power Systems}, 27(2):1138--1139, 2012.

\bibitem{JCC02}
R.~A. {Jabr}, A.~H. {Coonick}, and B.~J. {Cory}.
\newblock A primal-dual interior point method for optimal power flow
  dispatching.
\newblock {\em IEEE Transactions on Power Systems}, 17(3):654--662, 2002.

\bibitem{JMPG15}
C.~{Josz}, J.~{Maeght}, P.~{Panciatici}, and J.~C. {Gilbert}.
\newblock Application of the moment-sos approach to global optimization of the
  opf problem.
\newblock {\em IEEE Transactions on Power Systems}, 30(1):463--470, 2015.

\bibitem{KDS16}
B.~Kocuk, S.~S. Dey, and X.~A. Sun.
\newblock Strong socp relaxations for the optimal power flow problem.
\newblock {\em Oper. Res.}, 64(6):1177–1196, December 2016.

\bibitem{Las01}
J.B. Lasserre.
\newblock Global optimization with polynomials and the problem of moments.
\newblock {\em SIAM Journal on Optimization}, 11(3):796--817, 2001.

\bibitem{LavLow12}
J.~{Lavaei} and S.~H. {Low}.
\newblock Zero duality gap in optimal power flow problem.
\newblock {\em IEEE Transactions on Power Systems}, 27(1):92--107, 2012.

\bibitem{MSL15}
R.~Madani, S.~Sojoudi, and J.~Lavaei.
\newblock Convex relaxation for optimal power flow problem: Mesh networks.
\newblock {\em IEEE Transactions on Power Systems}, 30:199--211, 2015.

\bibitem{Cor76}
G.P. McCormick.
\newblock Computability of global solutions to factorable non-convex programs:
  Part i - convex underestimating problems.
\newblock {\em Mathematical Programming}, 10(1):147--175, 1976.

\bibitem{MolHis15}
D.~K. Molzahn and I.~A. Hiskens.
\newblock {Sparsity-Exploiting Moment-Based Relaxations of the Optimal Power
  Flow Problem}.
\newblock {\em IEEE Transactions on Power Systems}, 30(6):3168--3180, November
  2015.

\bibitem{MHLD13}
D.~K. {Molzahn}, J.~T. {Holzer}, B.~C. {Lesieutre}, and C.~L. {DeMarco}.
\newblock Implementation of a large-scale optimal power flow solver based on
  semidefinite programming.
\newblock {\em IEEE Transactions on Power Systems}, 28(4):3987--3998, 2013.

\bibitem{Par03}
P.~A. Parrilo.
\newblock Semidefinite programming relaxations for semialgebraic problems.
\newblock {\em Mathematical Programming Ser B.}, 96(2):293–320, 2003.

\bibitem{PGLIB}
{PGL}ib {O}ptimal {P}ower~{F}low {B}enchmarks.
\newblock {\em The {IEEE PES T}ask {F}orce on {B}enchmarks for {V}alidation of
  {E}merging {P}ower {S}ystem {A}lgorithms}, "accessed 02/2021".

\bibitem{baron}
N.V. Sahinidis and M.~Tawarmalani.
\newblock Baron 9.0.4: Global optimization of mixed-integer nonlinear programs.
\newblock {\em User's Manual}, 2010.

\bibitem{Torres1998AnIM}
G.~L. Torres and V.~Quintana.
\newblock An interior-point method for nonlinear optimal power flow using
  voltage rectangular coordinates.
\newblock {\em IEEE Transactions on Power Systems}, 13:1211--1218, 1998.

\bibitem{ipopt}
A.~W{\"a}chter and L.T. Biegler.
\newblock On the implementation of an interior-point filter line-search
  algorithm for large-scale nonlinear programming.
\newblock {\em Mathematical Programming}, 106(1):25--57, Mar 2006.

\bibitem{WMZT07}
H.~Wang, C.E. Murillo-S\'anchez, R.D. Zimmerman, and R.J. Thomas.
\newblock On computational issues of market-based optimal power flow.
\newblock {\em IEEE Transactions on Power Systems}, 22(3):1185--1193, 2007.

\bibitem{Wu1993ADN}
Y.C. Wu, A.~Debs, and R.~Marsten.
\newblock A direct nonlinear predictor-corrector primal-dual interior point
  algorithm for optimal power flows.
\newblock {\em IEEE Transactions on Power Systems}, 9:876--883, 1993.

\end{thebibliography}
     \end{document}